\documentclass[11pt,a4paper]{amsart}

\usepackage{amsmath,mathrsfs,amssymb}

\newtheorem{theorem}{Theorem}[section]
\newtheorem{lemma}[theorem]{Lemma}
\newtheorem{proposition}[theorem]{Proposition}
\newtheorem{corollary}[theorem]{Corollary}

\newtheorem{question}[theorem]{Question}

\theoremstyle{definition}

\theoremstyle{plain}

\newcommand{\GreenL}{\mathscr{L}}
\newcommand{\GreenR}{\mathscr{R}}
\newcommand{\GreenH}{\mathscr{H}}
\newcommand{\GreenD}{\mathscr{D}}

\newcommand{\BN}{\mathbb{N}}

\newcommand{\NO}{\textbf{NO}}
\newcommand{\YES}{\textbf{YES}}

\begin{document}
    \title[Uniform Decision Problems in Automatic Semigroups]{Uniform Decision Problems in\\ Automatic Semigroups}

\maketitle

\begin{center}
    Mark Kambites and Friedrich Otto\\

    \medskip

    Fachbereich Mathematik / Informatik, \   Universit\"at Kassel \\
    34109 Kassel, \  Germany \\

    \medskip

    \texttt{[kambites|otto]@theory.informatik.uni-kassel.de} \\

\bigskip

\end{center}

\begin{abstract}
We consider various decision problems for automatic semigroups, which
involve the provision of an automatic structure as part of the problem
instance. With mild restrictions on the automatic structure, which seem
to be necessary to make the problem well-defined, the uniform word
problem for semigroups described by automatic structures is decidable.
Under the same conditions, we show that one can also decide whether the
semigroup is completely simple or completely zero-simple; in the case
that it is, one can compute a Rees matrix representation for the semigroup,
in the form of a Rees matrix together with an automatic structure for its
maximal subgroup. On the other hand, we show that it is undecidable in
general whether a given element of a given automatic monoid has a right
inverse.
\end{abstract}

\section{Introduction}

Over the past two decades, one of the most successful and productive areas
of computational algebra has been the theory of \textit{automatic groups}.
Roughly speaking, the description of a group by an automatic structure
allows one efficiently to perform various computations involving the group,
which may be hard or impossible given only a presentation. Groups which
admit automatic structure also share a number of interesting structural
and geometric properties \cite{Epstein92}. More recently, many authors
have followed a suggestion of Hudson
\cite{Hudson96} by considering a natural generalisation to the
broader class of monoids or, even more generally, of
semigroups, and a coherent theory has begun to develop
\cite{Campbell99,Campbell01,Campbell00a,Descalco01,Duncan99,Hoffmann01,Hoffmann02,Hoffmann00b,Hoffmann01a,Otto99,Otto00,Otto00b,Otto01,Otto98}.

A number of authors have considered decision problems in automatic
semigroups; for example it has been established that the word problem for
an automatic semigroup is always decidable in quadratic time
\cite{Campbell01} and in certain cases is P-complete \cite{Lohrey04}.
In general, this research has assumed a fixed semigroup with automatic
structure, and asked what computations can be performed in the semigroup.
Since an automatic structure is a finite description of a (typically
infinite) semigroup, one is also able to consider decision problems in which
a semigroup, defined by means of an automatic structure, forms part of
the problem instance. Some such problems, such as the isomorphism problem
and the uniform word problem, are also studied by group theorists. Others
are more particular to semigroup theory; for example, one can ask if there
are algorithms to decide, given an automatic structure, whether the
semigroup described is a group, is completely simple, is completely-zero-simple,
has a left/right/two-sided zero or identity, is right/left/two-sided
cancellative and so forth.

The aim of this paper is to make a start upon addressing these issues.
We begin, in Section~\ref{sec_foundations} by developing a suitable
theoretical foundation for the study of uniform problems involving
automatic semigroups. In Sections~\ref{sec_unint} and~\ref{sec_int}
we present a number of algorithms for basic problems involving
automatic semigroups; these include, amongst others, the uniform word
problem and deciding the properties of right cancellativity, the
existence of an identity and the existence of a zero. In
Section~\ref{sec_rightinversesundec}, by contrast, we show that 
the existence of right inverses for a given element of a given
automatic semigroup is not, in general, decidable. Finally, in
Section~\ref{sec_compsimp}, we present algorithms to decide if a given
semigroup is completely simple or completely zero-simple and, in the
event that it is, obtain a Rees matrix decomposition of the semigroup.

While this paper is theoretical in nature, the research documented also
has a practical aspect. We aim not only to improve our understanding of
computational issues involving automatic semigroups, but also to develop
practical tools which will be of use to those working in the field. All
the algorithms documented in this paper, along with a number of others,
have been implemented by the first author using the GAP computer algebra
system \cite{GAP}. They will shortly be made available in the form of a
GAP package, as a resource for researchers in the area.

\section{Foundations}\label{sec_foundations}

In this section, we begin by recalling some basic definitions which we
shall require in the sections that follow. We then proceed to develop a
suitable formalism for the study of uniform problems involving automatic
structures. We assume a basic familiarity with finite automata; the
reader with no experience in this area is advised to consult a textbook
such as \cite{Lawson03}. Throughout this paper, we write $\BN_0$ and
$\BN_+$ to denote the sets of non-negative integers and strictly positive
integers respectively, and $\epsilon$ to denote the empty word.

\subsection{Synchronous automata recognising relations}

Let $A$ be a finite alphabet, let $\$ $ be a new symbol not in $A$, and let $A^\$ $ be the alphabet
$A \cup \lbrace \$ \rbrace$. We define a function
$\delta : A^* \times A^* \to (A^\$ \times A^\$)^*$ as follows. For
$m, n  \in \BN_0$ and $a_1, \dots a_m, b_1, \dots, b_n \in A$, let
$$\delta (a_1 \dots a_m, b_1 \dots b_n) = \begin{cases}
 (a_1, b_1) \dots (a_m, b_m) (\$, b_{m+1}) \dots (\$, b_n) &\text{ if } m < n \\
 (a_1, b_1) \dots (a_m, b_m) &\text{ if } m = n \\
 (a_1, b_1) \dots (a_n, b_n) (a_{n+1}, \$) \dots (a_m, \$) &\text{ if } m > n.
\end{cases} $$
Intuitively, $\delta$ rewrites \textit{pairs of words} over $A$ as
\textit{words of pairs} over $A^\$ $.

A \textit{synchronous automaton} over $A$ is a finite automaton over
$(A^\$ \times A^\$)^*$. We say that such an automaton \textit{recognises}
or \textit{accepts} a pair $(u, v) \in A^* \times A^*$ if it recognises
$\delta(u, v)$; it recognises a relation $R \subseteq A^* \times A^*$ if
it recognises exactly the image $\delta(R)$. A relation recognised by a
synchronous automaton is called $\textit{synchronously rational}$.
The following proposition summarises some basic properties of synchronously
rational relations, which we shall use throughout this paper without further
comment; all are proved in \cite{Epstein92}.
\begin{proposition}
The class of synchronously rational relations is closed under
intersection, union, complement, composition and inversion. Synchronously
rational relations are rational transductions, and hence preserve
rational languages; in particular the projection of a synchronously
rational relation onto either coordinate is a rational language, and the
image of an element under a synchronously rational relation is a rational
language.

Moreover, all of these operations are effectively computable.
\end{proposition}
We shall also make use of the following fact, again without further comment.
\begin{proposition}
There is an algorithm which, given two synchronous automata, decides if
they accept the same relation.
\end{proposition}
\begin{proof}
Two synchronous automata accept the same relation exactly if, when
considered as normal finite automata, they accept the same language
over $(A^\$ \times A^\$)^*$. The latter property is well-known to
be decidable (for example, by computing minimal deterministic automata).
\end{proof}

\subsection{Automatic structures and interpretations}

A \textit{(synchronous) pre-automatic structure} $\Gamma$ consists of
\begin{itemize}
\item[(i)] a finite set $A(\Gamma)$ of \textit{generators};
\item[(ii)] a finite automaton recognising a language $L(\Gamma)$ over $A(\Gamma)$;
\item [(iii)] a synchronous automaton recognising a relation $L_=(\Gamma)$ on
      $A(\Gamma)^*$, which is contained in $L(\Gamma) \times L(\Gamma)$; and
\item[(iv)] for each $a \in A(\Gamma)$, a synchronous automaton recognising a
relation $L_a(\Gamma)$ on $A(\Gamma)^*$ which is contained in $L(\Gamma) \times L(\Gamma)$.
\end{itemize}
Where only one automatic structure is under discussion, we shall for
brevity write simply $A$, $L$, $L_=$ and $L_a$ in place of
$A(\Gamma)$, $L(\Gamma)$, $L_=(\Gamma)$ and $L_a(\Gamma)$ respectively.

An \textit{interpretation} of a pre-automatic structure
with respect to a semigroup $S$ is a morphism $\sigma : A^* \to S$ such
that
\begin{itemize}
\item[(i)] $\sigma(L) = S$;
\item[(ii)] for $u, v \in L$ we have $(u,v) \in L_=$ if and only if $\sigma(u) = \sigma(v)$; and
\item[(iii)] for each $a \in A$ and $u, v \in L_a$ we have $(u,v) \in L_a$ if 
            and only if $\sigma(ua) = \sigma(v)$.
\end{itemize}
If such an interpretation exists, we say that the pre-automatic structure
is an \textit{automatic structure for $S$}, or that $S$ is
\textit{described} by the automatic structure. A semigroup
$S$ is called \textit{automatic} if it is described by some
automatic structure. If, in addition, the interpretation restricts to
a bijection from $L$ to $S$, we say that the automatic structure has
\textit{uniqueness}, or is an \textit{automatic cross-section} for $S$.

\begin{proposition}\label{prop_uniquesemigroup}
Two semigroups described by the same automatic structure
are necessarily isomorphic.
\end{proposition}
\begin{proof}
Suppose $\sigma : A^* \to S$ and $\tau : A^* \to T$ are interpretations
of the same automatic structure. Define a map $\rho : S \to T$
by setting $\rho(s)$ to be the unique $t \in T$ such that $t = \tau(w)$ for
some $w \in L$ with $\sigma(w) = s$. It is straightforward to verify that
this map is a well-defined isomorphism of the semigroups.
\end{proof}
Proposition~\ref{prop_uniquesemigroup} tells us that an
automatic structure uniquely defines a semigroup $S$ up to isomorphism,
but it does \textbf{not}
guarantee that it uniquely defines an interpretation, even up to
isomorphism. It is not clear whether every automatic structure contains
sufficient information to associate to each generator a member of
the language of representatives which represents the same element in
every interpretation.

We will need to perform computations not only with automatic
structures but also with interpretations; for this we require
a finite way of encoding of an interpretation.
Formally, we say that an \textit{assignment
of generators} for the automatic structure with respect to a semigroup $S$
is a function $\iota : A \to L$ with the property that there exists an
interpretion $\sigma : A^* \to S$ such that $\sigma(a) = \sigma(\iota(a))$
for all $a \in A$; such an interpretation is said to be \textit{consistent}
with $\iota$. Two assignments of generators $\iota$ and $\iota'$ are called
\textit{equivalent} if $(\iota(a), \iota'(a)) \in L_=$ for every
generator $a$.
An \textit{interpreted automatic structure} for a semigroup $S$ is a
pair $(\Gamma, \iota)$ of an automatic structure together with an
assignment of generators with respect to $S$. This terminology is
justified by the following proposition which is straightforward to prove.
\begin{proposition} \label{prop_inter_iff_assign}
An assignment of generators is consistent with a unique interpretation (up to
isomorphic permutation of the semigroup described). Moreover, equivalent
assignments of generators are consistent with the same interpretation.

Conversely, an interpretation for an automatic structure is
consistent with an assignment of generators, which is unique up to equivalence.
\end{proposition}

Given an automatic structure and a word $w = a_1 a_2 \dots a_n \in A^+$ with
each $a_i \in A$, we define
$$L_w = L_{a_1} \circ L_{a_2} \circ \dots \circ L_{a_n}$$
where $\circ$ denotes composition of relations.
We extend this definition to the whole of $A^*$ by letting $L_\epsilon = L_=$.
By our observations above, we can compute a synchronous automaton
recognising the language $L_w$ for any $w \in A^*$.
It is readily verified that for any interpretation
$\sigma : A^* \to S$ of the automatic structure, we have
$$L_w = \lbrace (u, v) \in L \mid \sigma(uw) = \sigma(v) \rbrace.$$

\subsection{Interpreted vs uninterpreted automatic structures}

The distinction between automatic structures with and without interpretation
is important for two reasons. Firstly, certain problems, such
as the uniform word problem, involve elements of the semigroup expressed 
as words in the generators as part of the problem instance; these problems
are not necessarily well-defined in the absence of an interpretation.
Secondly, even invariant properties
of the semigroup, such as whether it has an identity, can be more straightforward
to test when an assignment of generators is provided. We presently lack
an algorithm which, given an automatic structure,
computes an assignment
of generators; in general it is not clear whether such an algorithm exists.
\begin{question}
Is there an algorithm which, given an automatic structure, finds an
assignment of generators which is consistent with some interpretation?
\end{question}
However, there is a very large class of semigroups
for which automatic structures admit essentially unique interpretations, and
in which such an algorithm does exist.

We say that two elements $s$ and $t$ of a semigroup $S$ are \textit{right
translationally equivalent} if $xs = xt$ for every element $x \in S$. Recall
that a semigroup is called \textit{left reductive} if no two distinct elements
are right translationally equivalent \cite[p.~84]{Petrich73}. The class of
left reductive semigroups is very large, including all monoids, left
cancellative semigroups, inverse semigroups and of course groups. An example
of a semigroup which is \textbf{not} left reductive is a non-trivial semigroup
all of whose elements are left zeros.

\begin{proposition}\label{prop_righttransequiv}
Given an automatic structure, right translational equivalence is independent
of the choice of interpretation, that is, words $u$ and $v$ represent right
translationally equivalent elements in all interpretations or in none.

Moreover, there is an algorithm which, given
as input an automatic structure and two words $u$ and $v$ in the generators,
decides if $u$ and $v$ represent right translationally equivalent elements.
\end{proposition}
\begin{proof}
It is readily verified that $u$ and $v$ are right translationally equivalent
in any interpretation if and only if the languages $L_u$ and $L_v$ are equal.
These languages can be computed independent of the interpretation, by
composition of languages of the form $L_a$ for various $a \in A$. We can
then solve the problem by testing them for equality.
\end{proof}
\begin{corollary}\label{cor_leftredinter}
An automatic structure for a left reductive semigroup admits a unique
interpretation (up to isomorphism of the semigroup described). Moreover,
there is an algorithm which, given as input an automatic structure for
a left reductive semigroup, computes an assignment of generators.
\end{corollary}
\begin{proof}
Suppose $\iota$ and $\iota'$ are assignments of generators and $a \in A$.
Then $\iota(a)$ is right translationally equivalent to $a$, and hence to
$\iota'(a)$. Since the semigroup is left reductive, it follows that
$\iota(a)$ and $\iota'(a)$ represent the same element, and so we must have
$(\iota(a), \iota'(a)) \in L_=$. Thus, $\iota$ and $\iota'$ are equivalent,
and so by Proposition~\ref{prop_inter_iff_assign}, they are consistent with
the same unique interpretation.

To compute the interpretation, for each generator $a$ we enumerate all words
in $L$ until we find a word $w \in L$ which is right translationally
equivalent to $a$, and set $\iota(a) = w$.
\end{proof}

\subsection{Modifying automatic structures}

We shall make use of the following propositions, which allow us, given an
interpreted automatic structure, to obtain automatic structures with
nicer properties, for the same semigroup. They are essentially algorithmic
restatements of results from \cite{Campbell01} and \cite{Epstein92}.

\begin{proposition}\label{prop_modify}
There is an algorithm to solve the following problem:

\textbf{Instance:} an interpreted automatic structure with language of
representatives $L$ admitting an interpretation $\sigma : A^* \to S$ and
a finite automaton recognising a regular language $K \subseteq A^*$ such
$K \setminus L$ is finite and $\sigma(K) = S$;

\textbf{Problem:} compute an automatic structure admitting the same
interpretation and with language of representatives $K$.
\end{proposition}
\begin{proof}
By \cite[Propositions~5.3 and 5.7]{Campbell01} there exists such an
automatic structure. Moreover, the proofs of those results give an
effective method of construction.
\end{proof}

\begin{proposition}\label{prop_uniqueness}
There is an algorithm which, given as input an interpreted automatic
structure admitting an interpretation $\sigma : A^* \to S$ which is
injective when restricted to the generators, computes an automatic
structure such that
\begin{itemize}
\item[(i)] the new automatic structure admits the same interpretation $\sigma$;
\item[(ii)] the restriction of $\sigma$ to the language of representatives is bijective; and
\item[(iii)] the language of representatives contains the generators.
\end{itemize}
\end{proposition}
\begin{proof}
By Proposition~\ref{prop_modify} we may obtain an automatic structure
whose language of representatives $L$ contains the generators. In the
event that the identity is a generator, we may also assume that $L$ does
not contain the empty word.

Using the argument from \cite[Theorem~2.5.1]{Epstein92}, we can construct
an automaton recognising the language $L'$ of words in $L$ which are
minimal with respect to the shortlex order (see \cite[Section~2.5]{Epstein92})
amongst words in $L$ representating the same element
of the semigroup. Clearly, the generators will be contained in $L'$. Now by
Proposition~\ref{prop_modify}, we obtain a suitable automatic structure with
language of representatives $L'$.
\end{proof}

The previous proposition combines with the next one, to give us a more
concise way to encode an interpreted automatic structure.
\begin{proposition}\label{prop_shortencoding}
An automatic structure with the generators in the language of representatives admits a unique
interpretation up to equivalence, and there is an algorithm which, given such an
automatic structure, computes the interpretation.
\end{proposition}
\begin{proof}
The identity function on the set of generators serves as an assignment
of generators. It is easy to see that any other assignment of generators
must be equivalent to this one.
\end{proof}
In view of Propositions~\ref{prop_uniqueness} and \ref{prop_shortencoding},
we may describe an interpretation of an automatic structure by including
generators in the language of representatives; this avoids the need for
explicit reference to the assignment of generators.

\section{Algorithms for Uninterpreted Automatic Structures}\label{sec_unint}

In this section, we consider decision problems for which the instance is
an automatic structure without any further information about the assignment
of generators.

First, we recall from \cite[Section~5.1]{Epstein92} that
one can decide, starting only from a collection of synchronous automata,
whether they form an automatic structure for a group. In particular, one
can verify whether a given semigroup automatic structure describes a
group. It is interesting to note that the corresponding
problem with input specified as a finite presentation is known to be undecidable
\cite{Narendran91}. One wonders whether this difference results from
(i) the difference in method of presentation, or (ii) the smaller class
of semigroups under consideration. In particular, one might ask the
following question.
\begin{question}
Is there an algorithm to solve the following problem?

\textbf{Instance:} a finite presentation for a semigroup which is known to
                   be automatic;

\textbf{Problem:} decide if the semigroup is a group.
\end{question}
It is not presently known whether there is an algorithm which, given a
finite presentation for an automatic semigroup, computes
an automatic structure for the semigroup. If, in fact,
there is such an algorithm, then the answer to the previous question is
necessarily positive.

Another property which is easily decided is that of right cancellability.
\begin{proposition}
Let $\Gamma$ be an automatic cross-section, and $\sigma : A^* \to S$ an
interpretation. Let $w \in A^*$.
Then $\sigma(w)$ is right cancellable
in $S$ if and only if $L_w \circ L_w^{-1}$ is the diagonal relation on $L$.
\end{proposition}
\begin{proof}
Let $\sigma : A^* \to S$ be an interpretation. If
$L_w \circ L_w^{-1}$ is \textbf{not} diagonal then there exist distinct\
words $u, v \in L$ such that $(u,v) \in L_w \circ L_w^{-1}$,
that is, such that there exists $x \in L$ with $(u, x), (v,x) \in L_w$.
But now $\sigma(u) \sigma(w) = \sigma(x) = \sigma(v) \sigma(w)$. But
since the automatic structure is a cross-section, we have that $\sigma(u) \neq \sigma(v)$,
so that $\sigma(w)$ is not right cancellable.

Conversely, if $\sigma(w)$ is not right cancellable, then $a \sigma(w) = b \sigma(w)$
for some distinct elements $a, b \in S$. Let $u$, $v$ and $x$ be words in $L$
representing $a$, $b$ and $a \sigma(w)$ respectively. Then we have
$(u, x), (v,x) \in L_w$, so that $(u, v) \in L_w \circ L_w^{-1}$. Clearly
$u \neq v$, so $L_w \circ L_w^{-1}$ is not diagonal.
\end{proof}
\begin{corollary}\label{cor_rightcanc_decidable}
There is an algorithm for the following problem:

\textbf{Instance:} an automatic structure for a semigroup;

\textbf{Problem:} decide whether the semigroup is right cancellative.
\end{corollary}
\begin{proof}
It suffices to check that each generator $a \in A$ is right cancellable
by computing $L_a \circ L_a^{-1}$ and comparing with the diagonal relation
on $L$.
\end{proof}

\begin{proposition}\label{prop_leftzero_decidable}
There is an algorithm for the following problem:

\textbf{Instance:} an uninterpreted automatic structure;

\textbf{Problem:} decide whether the semigroup described has a
                  left zero, and if so find a finite automaton
                  representing the language of all words in $L$
                  representing left zeros.
\end{proposition}
\begin{proof}
For a word $z \in L$ we have that $z$ is a left
zero if and only if $za = z$ in the semigroup for all generators $a \in A$,
that is, if $(z,z) \in L_a$ for all generators $a \in A$. Hence, 
computing the intersection of all the languages $L_a$ with the diagonal
relation, and taking the projection onto one coordinate, gives the
language of all words in $L$ representing left zeros. If this language is
empty, output \NO. If not, output \YES\ and the automaton computed for
the language.
\end{proof}

\begin{proposition}\label{prop_zero_decidable}
There is an algorithm for the following problem:

\textbf{Instance:} an uninterpreted automatic structure;

\textbf{Problem:} decide whether the semigroup described has a
                  zero, and if so find the word in $L$ representing it.
\end{proposition}
\begin{proof}
By Proposition~\ref{prop_uniqueness}, we may assume that the automatic
structure has uniqueness. First compute the language of left zeros, as
above. Since a zero in a semigroup must be the unique left zero, check
that this language contains exactly one element; if not output \NO.
Otherwise, check that $z$ represents a right zero by checking that $uz = z$ for
all words $u \in L$, that is, that $L_z = L \times \lbrace z \rbrace$.
If so, output \YES, otherwise output \NO.
\end{proof}

\begin{question}
Is there an algorithm for the following problem?

\textbf{Instance:} an uninterpreted automatic cross-section $(A, L)$;

\textbf{Problem:} decide if the semigroup presented is a monoid.
\end{question}

\section{Basic Algorithms for Interpreted Automatic Structures}\label{sec_int}

In this section we consider various basic algorithmic problems which
start with an interpreted automatic structure. We begin with the
following simple proposition.

\begin{proposition}\label{prop_findwordinl_decidable}
There is an algorithm to solve the following problem:

\textbf{Instance:} an interpreted automatic structure and a non-empty word $u$ over
                   the generators;

\textbf{Problem:} find a word in the language of representatives, which
                  represents the same element as $u$.
\end{proposition}
\begin{proof}
We use induction on the length of $u$. The base case is given by the
assignment of generators. Now in general, write $u = va$ where $a$ is
a generator. By the inductive hypothesis, we can find a word $w$ in
the language of representatives representing the same element as $v$.
Now a word $x$ in the language of representatives represents $u$ if and
only if $(w, x) \in L_a$.
Using standard operations on finite automata, it is straightforward
to find such an $x$.
\end{proof}
\begin{corollary}[Uniform Word Problem]
There is an algorithm to solve the following problem:

\textbf{Instance:} an interpreted automatic structure and two words $u$
                   and $v$ over the generators;

\textbf{Problem:} decide if $u$ and $v$ represent the same element.
\end{corollary}
\begin{proof}
By Proposition~\ref{prop_findwordinl_decidable}, we can compute elements $w_u$ and
$w_v$ in $L$ representing $u$ and $v$ respectively, and then check whether
$(w_u, w_v) \in L_=$.
\end{proof}

With an interpretation available, it becomes an easy task to decide
whether an automatic semigroup has an identity, and if so to locate
a representative for it.

\begin{proposition}\label{prop_identity_decidable}
There is an algorithm for the following problem:

\textbf{Instance:} an interpreted automatic structure for a semigroup;

\textbf{Problem: } Decide whether the semigroup described is a monoid,
                   and if so, find a word in L representing the identity.
\end{proposition}
\begin{proof}
By Proposition~\ref{prop_uniqueness} we may assume that the automatic
structure has uniqueness and that the language of representatives contains
the generators.

Note that for $w \in L$ and $a \in A$, we have $wa=a$ in the semigroup
if and only if $(w, a) \in L_a$. It follows that we can easily compute
the set of words in $L$ which stabilise the generators on the left, that
is, which represent left identities in the monoid. Call this set $K$.

Since the automatic structure is a cross-section, every word in $K$ represents a different
left identity of the semigroup. A monoid can have only one left identity,
so if $K$ is not a singleton set, output \NO. Otherwise let $e$ be the
unique word in K. Now check if $e$ is a right identity, by verifying that
for each generator $a \in A$ we have $ae = a$ (or alternatively that
$L_e = L_=$). If so, output \YES\ and the word e, otherwise output \NO.
\end{proof}

Recall that a \textit{right inverse} [\textit{left inverse}] of a
monoid element $s$ is a monoid element $t$ such that $st = 1$
[respectively, $ts = 1$]. An element with a left and a right
inverse (which necessary coincide) is called a \textit{unit}; the
set of all units in a monoid forms a subgroup of the monoid.

\begin{proposition}\label{prop_leftinverse_decidable}
There is an algorithm to solve the following problem:

\textbf{Instance:} an interpreted automatic structure
describing a monoid, and a word $w$ in the generators;

\textbf{Problem:} decide whether $w$ has a left inverse in the monoid,
and obtain an automaton recognising the language of representatives in
$L$ of its left inverses.
\end{proposition}
\begin{proof}
By Proposition~\ref{prop_identity_decidable}, we can compute a word
$e \in L$ representing the identity. Now a left inverse of $w$ is an
element represented by a word $w' \in L$ such that $(w', e) \in L_w$.
It is straightforward to compute the language of such and check whether
it is empty.
\end{proof}

In contrast, we shall see in the Section~\ref{sec_rightinversesundec}
below that it is not in general possible to decide whether a word $w$
represents an element with a right inverse.

\begin{proposition}\label{prop_unit_decidable}
There is an algorithm for the following problem:

\textbf{Instance:} an interpreted automatic structure describing a monoid
and a word $w$ in the generators;

\textbf{Problem:} decide whether $w$ represents a unit in the
monoid.
\end{proposition}
\begin{proof}
By Proposition~\ref{prop_uniqueness}, we may assume that the automatic
structure has uniqueness.
First use the method above to check that $w$ has a left inverse, and obtain
the language of its left inverses. If $w$ is to be a unit then it can
have only one left inverse, so if this language has more than one element,
output \NO. Otherwise, let $w'$ be the unique element of the language,
and check that $w w' = e$.
\end{proof}

Let $S$ be a semigroup, and $0$ a new symbol not in $S$. We define a
new semigroup $S^0$ with set of elements $S \cup \lbrace 0 \rbrace$,
and multiplication given by
$$xy = \begin{cases} 0 &\text{ if } x = 0 \text{ or } y = 0 \\
                     \text{the $S$-product } xy &\text{ otherwise}.\end{cases}$$
The semigroup $S^0$ is called \textit{$S$ with an adjoined zero}. The
following result is essentially an algorithmic restatement of
\cite[Proposition~3.13]{Campbell01}.

\begin{proposition}\label{prop_adjoinzero}
There is an algorithm for the following problem:

\textbf{Instance:} an interpreted automatic structure $\Gamma$
                   describing a semigroup $S$;

\textbf{Problem:} compute an automatic structure for the semigroup $S^0$.
\end{proposition}

We saw above that one can decide, given even an uninterpreted automatic
structure, whether the semigroup described is right cancellative. However,
the following question remains open.

\begin{question}
Is there an algorithm to decide, from an interpreted
automatic structure, whether the semigroup is left cancellative
or cancellative?
\end{question}

\section{Undecidability of the existence of a right inverse.}\label{sec_rightinversesundec}

The aim of this section is to demonstrate the existence of automatic
monoids for which there is no algorithm to decide, given a word $u$ over
the generating set, whether $u$ represents an element with a right
inverse in the monoid. We use without further comment a number of standard results
from the theory of string-rewriting; these can be found in \cite{Book93}.

Let $M = (Q,\Sigma,B,q_0,q_a,\delta)$ be a deterministic single-tape Turing
machine,
where $Q$ is the finite set of states,
$\Sigma$ is the finite tape alphabet,
$B\notin \Sigma$ is the blank symbol,
$q_0\in Q$ is the initial state,
$q_a\in Q$ is the final (accepting) state and
$$\delta : (Q \times \Sigma) \to (Q \times \Sigma \times \lbrace \lambda, \rho \rbrace)$$
is the transition function.
We assume that $M$ halts if and when it enters the state $q_a$,
and that the tape of $M$ is unrestricted only to the right.
Given a word $w\in\Sigma^*$ as input,
the corresponding initial configuration of $M$ is $q_0w$,
where we assume that $w$ occupies the prefix of length $|w|$ of the tape.
The word $w$ is accepted by $M$ if and only if the computation of $M$
that starts from this initial configuration finally ends in the final state.
As the tape is limited on the left, during a computation the head of $M$ cannot ever move
to the left of its start position.
By $L(M)$ we denote the set of all words that are accepted by~$M$.

From $M$ we now construct a finite string-rewriting system $R_M$ on the
alphabet
$$\Gamma := Q\cup\Sigma\cup\overline{\Sigma}\cup\{d,h,\bar{h}\},$$
where $\overline{\Sigma} := \{\,\bar{s}\mid s\in\Sigma\,\}$ is a marked
copy of~$\Sigma$, and $d$, $h$, and $\bar{h}$ are three new symbols.
The system $R_M$ consists of the following rules:

$$
\begin{array}{llclllll}
{\rm (1)}\; & qad & \to & \bar{b}p & \mbox{if}\quad & \delta(q,a) & =\quad & (p,b,\rho),\\
{\rm (2)} & qhd & \to & \bar{b}ph & \mbox{if} & \delta(q,B) & = & (p,b,\rho),\\
{\rm (3)} & \bar{c}qad & \to & pcb & \mbox{if} & \delta(q,a) & = & (p,b,\lambda),\;c\in\Sigma,\\
{\rm (4)} & \bar{c}qhd & \to & pcbh\quad & \mbox{if} & \delta(q,B) & = & (p,b,\lambda),\;c\in\Sigma,\\
{\rm (5)} & q_aad  & \to & q_a & \mbox{for} & a\in\Sigma,\\
{\rm (6)} & \bar{a}q_ahd  & \to & q_ah & \mbox{for} & a\in\Sigma,\\
{\rm (7)} & \bar{h}q_ahd & \to & \varepsilon,\\
{\rm (8)} & abd  & \to & adb &  \mbox{for} & a,b\in\Sigma,\\
{\rm (9)} & ahd & \to & adh & \mbox{for} & a\in\Sigma.\\
\end{array}
$$

We define a binary relation $>$ on $\Gamma^*$ as follows:
$$\begin{array}{lcl}
u > v \quad & \iff \quad & |u|_d > |v|_d\quad\mbox{or}\quad|u|_d = n = |v|_d,\\
 & & u = u_0du_1d\dots u_{n-1}du_n,\;v=v_0dv_1d\dots v_{n-1}dv_n,\;\mbox{and}\\
 & & \exists\,j\,: |u_j| > |v_j|\;\mbox{and}\;|u_{i}|=|v_i|\;
\mbox{for all}\;0\le i < j.
\end{array}$$
It is easily seen that $>$ is the strict part of a
partial ordering on~$\Gamma^*$ that is well-founded.
Further, whenever $u \rightarrow_{R_M} v$ holds, then $u > v$.
Thus, $R_M$ does not generate any infinite reduction sequences,
that is, $R_M$ is noetherian.
As there are no non-trivial overlaps between the left-hand sides of the rules of~$R_M$
(recall that $M$ is assumed to be a deterministic Turing machine),
we see that $R_M$ is also confluent.
Thus, $R_M$ is a finite convergent system,
which implies that the set ${\sf IRR}(R_M)$ of irreducible words mod~$R_M$
is a regular set of normal forms for the monoid $S_M$ that is given through the
finite presentation $(\Gamma;R_M)$.
\vspace{+0.2cm}

\noindent
{\bf Claim~1.} The language ${\sf IRR}(R_M)$ forms the language of
               representatives for an automatic structure for $S_M$.

\noindent
\begin{proof}
It remains to show that, for each symbol $a\in\Gamma$,
the right-multiplication relation $L_a$ is synchronously regular.
For each rule $(\ell\to r)\in R_M$, we see that $\ell$ ends with the symbol~$d$.
Thus, for each symbol $a\in\Gamma\setminus\{d\}$,
if $u\in{\sf IRR}(R_M)$, then also $ua\in{\sf IRR}(R_M)$.
Thus, for all these letters the corresponding language $L_a$ is clearly
synchronously regular.

It remains to consider the language~$L_d$.
Let $u\in{\sf IRR}(R_M)$.
Then there are three mutually exclusive cases:
\begin{enumerate}
\item $ud$ is also irreducible modulo $R_M$;
\item $u = xaz$, where $a \in \Sigma$, $z\in(\Sigma\cup\{h\})^+$,
and $xad\in{\sf IRR}(R_M)$, which implies that $xadz$ is the irreducible
descendant of $ud = xazd$ modulo $R_M$.
\item $u = xyaz$ for some factors $x,y,z\in{\sf IRR}(R_M)$ and a
letter $a\in\Sigma\cup\{h\}$
satisfying the following conditions:
\begin{enumerate}
\item $yad$ is the left-hand side of one of the rules of type (1) to (7) of $R_M$,
\item if $a\in\Sigma$, then $z\in (\Sigma\cup\{h\})^*$,
and if $a=h$, then $z=\varepsilon$.
\end{enumerate}
In this case the irreducible descendant of $ud = xyazd$ is the word $xrz$,
where $r$ is the right-hand side of the rule with left-hand side $yad$.
\end{enumerate}
In the first case $(u, ud) \in L_d$, in the second case $(xaz,xadz) \in L_d$,
while in the third case $(xyaz,xrz) \in L_d$. It is easily verified that these
observations imply that the language $\delta(L_d)$ is regular so
that $L_d$ is synchronously regular as required.
\end{proof}

Thus, the monoid $S_M$ is automatic.
\vspace{+0.2cm}

\noindent
{\bf Claim~2.} For each word $w\in\Sigma^*$, $w \in L(M)$ if and only if
there exists an integer $n>0$ such that
$$\bar{h}q_0wh\cdot d^n\rightarrow^+_{R_M}\varepsilon.$$

\begin{proof}
The reductions modulo $R_M$ essentially just simulate the steps of the Turing machine~$M$.
Each step of the simulation digests one occurrence of the symbol~$d$.
Further occurrences of the symbol $d$ are needed to reduce the encoding $\bar{h}\bar{x}q_ayh$
of the final configuration of $M$ to the word $\bar{h}q_ah$, and another occurrence of the symbol
$d$ is needed for the final step reducing $\bar{h}q_ah$ to the empty word.
\end{proof}

\vspace{+0.2cm}

\noindent
{\bf Claim~3.} For each word $w\in\Sigma^*$, $\bar{h}q_0wh$ is right-invertible in $S_M$
if and only if $w\in L(M)$.

\begin{proof}
If $w\in L(M)$, then $\bar{h}q_0wh$ is right-invertible in $S_M$ by Claim~2.
Conversely, if $\bar{h}q_0wh$ is right-invertible in $S_M$,
then there exists an irreducible word $z$ such that $\bar{h}q_0wh\cdot z\rightarrow^+_{R_M} \varepsilon$ holds.
As $\bar{h}q_0wh$ and $z$ are both irreducible,
rewrite steps can only be applied across the border of these two factors.
It follows that $z = d^n$ for some positive integer~$n$,
which in turn implies by Claim~2 that $w\in L(M)$ holds.
\end{proof}

Thus, the halting problem for the Turing machine~$M$ reduces to the
problem of finding a right inverse for an element in the automatic
monoid~$S_M$, giving the following result.

\begin{theorem}
There exists a finitely presented automatic monoid $S$ and an
automatic structure for $S$ for which the following problem is
in general undecidable:

\textbf{Instance:} a word $w$ in the language of representatives;

\textbf{Problem:} decide if $w$ represents an element with a right inverse
in $S$.
\end{theorem}

\section{Completely simple and completely zero-simple semigroups.}\label{sec_compsimp}

In this section, we present algorithms to decide if a given automatic structure
represents a completely simple or completely zero-simple semigroup and, in
the case that it does, to compute a Rees matrix decomposition for the semigroup.
As well as being of interest in its own right, this shows that the use of
automatic structures to describe semigroups facilitates the decision of
quite complex structural properties.

Recall that a \textit{primitive idempotent} in a semigroup $S$ is an
idempotent $e$ with the property that for any non-zero idempotent $f$
such that $ef = fe = f$, we have $e = f$. A semigroup is
called \textit{completely simple} if it has a primitive idempotent and no proper
ideals. A semigroup with zero is called \textit{completely zero-simple}
if the zero is the only proper ideal, and it has no infinite descending chains
of idempotents. For a detailed introduction to the theory of completely
simple and completely zero-simple semigroups, see \cite[Chapter~3]{Howie95}.

The following construction, due to Rees \cite{Rees40}, provides a way
to describe completely zero-simple and completely simple semigroups.
Let $G$ be a group and $I$ and $\Lambda$ be sets. Let $0$ be a symbol
not in $G$, and let $P$ be a $\Lambda \times I$ matrix with entries drawn
from $G \cup \lbrace 0 \rbrace$. The \textit{Rees matrix semigroup with zero}
$M^0(G; I, \Lambda; P)$ is the semigroup with set of elements
$$\left( I \times G \times \Lambda \right) \cup \lbrace 0 \rbrace$$
and multiplication given by
$$(i, g, \lambda) (j, h, \mu) =
 \begin{cases} (i, g P_{\lambda j} h, \mu) &\text{ if } P_{\lambda j} \in G; \\
               0 &\text{ if } P_{\lambda j} = 0
 \end{cases}$$
and $x0 = 0 = 0x$ for all elements $x$.

The matrix $P$ is called \textit{regular} if every row and every column
contains a non-zero entry. If $P$ contains no zero entries at all, then
the set of non-zero elements of $M^0(G; I, \Lambda; P)$ forms a subsemigroup,
called a \textit{Rees matrix semigroup (without zero)} and denoted $M(G; I, \Lambda; P)$.

The following theorem is usually attributed to Rees, although it was essentially
prefigured by Suschkewitz \cite{Suschkewitz28}. It provides the connection
between completely simple semigroups and Rees matrix constructions.

\begin{theorem}[Suschkewitz 1928, Rees 1940]\label{thm_rees}
Let $G$ be a group, $I$ and $\Lambda$ sets, and $P$ a regular
$\Lambda \times I$ matrix over $G \cup \lbrace 0 \rbrace$. Then
$M^0(G; I, \Lambda; P)$ is a completely zero-simple semigroup.
If $P$ contains no zero entries, then $M(G; I, \Lambda; P)$ is a
completely simple semigroup.

Conversely, every completely simple or completely zero-simple semigroup
is isomorphic to one of this form.
\end{theorem}

Completely simple and completely zero-simple semigroups and 
Rees matrix constructions are of fundamental importance in the theory of semigroups.
Various authors have considered the relationship between Rees matrix
constructions and automaticity properties. In \cite{Campbell00a}, it is
shown that a finitely generated completely simple semigroup is automatic
if and only if its maximal subgroups are automatic; a consequence of
results in \cite{Descalco01} is that the same applies in the completely
zero-simple case. However, the methods used in these papers are essentially
non-constructive, in that they presuppose prior knowledge of the Rees
matrix representation for the semigroup.

In this section, we present two algorithms relating automatic structures
to completely simple and completely zero-simple semigroups; the first
takes as input an interpreted automatic structure, and decides whether
the semigroup represented is completely zero-simple. Our second algorithm takes as
input an automatic structure presupposed to represent a completely zero-simple
semigroup, and computes a Rees matrix representation for the semigroup;
the latter takes the form of an automatic structure for the (necessarily
unique up to isomorphism) maximal subgroup, and a sandwich matrix of words
from the language of representatives. We also show how to apply these
results to completely simple semigroups without zero.

\subsection{Deciding complete simplicity and complete zero-simplicity}

In this section, we show that there is an algorithm which, given an
automatic structure, decides whether the semigroup described is
completely simple or completely zero-simple. We shall require the following
lemma.
\begin{lemma}\label{lemma_inverse_conditions}
There is an algorithm for the following problem:

\textbf{Instance:} an interpreted automatic structure representing a
  semigroup $S$, and two words $w, e$ in the language of representives
  such that $e$ represents an idempotent in $S$;

\textbf{Problem:} decide which of the following comprehensive and
mutually exclusive conditions applies:
\begin{itemize}
\item[(A)] $w$ has an infinite number of left inverses which respect to $e$.
\item[(B)] $w$ has a finite number of left inverses with respect to $e$,
            one of which is also a right inverse;
\item[(C)] $w$ has a finite number of left inverses with respect to $e$,
             none of which is a right inverse.
\end{itemize}
\end{lemma}
\begin{proof}
By Proposition~\ref{prop_uniqueness}, we may assume that the automatic
structure has uniqueness. We begin by computing the language $K$ of left inverses for $w$ with 
respect to $e$. This is the set of all words $w'$ such that $(w', e) \in 
L_w$. If $K$ is infinite, output (A). Otherwise, check each $w' \in K$
in turn to see if $w w' = e$ in the semigroup. If one does, output 
$(B)$; otherwise, output $(C)$.
\end{proof}

\begin{theorem}\label{thm_compzero_decide}
There is an algorithm for the following problem:

\textbf{Instance:} an interpreted automatic structure for a semigroup $S$;

\textbf{Problem:} decide if $S$ is completely zero-simple.
\end{theorem}

\begin{proof}
    By Proposition~\ref{prop_uniqueness}, we may assume that we are
    given an interpreted automatic structure with uniqueness, and that
the language of representatives contains the generators. By
    Proposition~\ref{prop_zero_decidable}, we can check that the semigroup
    has a zero. If not, output \NO. If so, let $z$ be the (unique) word
    in $L$ representing the zero. Check which, if any, of the generators
    represent zero.

    Now for every generator $a$ which does not represent zero, calculate the set of words in $L$
    representing elements which stabilise $a$ on the left (that
    is, the projection onto
    the first coordinate of $L_a \cap (A^* \times \lbrace a \rbrace)$). Call it $SL_a$. In
a completely zero-simple semigroup, the elements which stabilise a non-zero element
$s$ on the left are exactly the idempotents in the $\GreenR$-class of $s$; it
follows that (i) there must be one such, (ii) there can only be finitely
many such (by the Main Theorem of \cite{Ayik99}) and (iii) they are all
idempotents. Check these three conditions, and if any fail, output \NO.

    Let $E$ be the union of the $SL_a$'s; thus, $E$ is a (finite) set
    of non-zero idempotents. (If
    $S$ is indeed completely zero-simple then it follows easily from
    Theorem~\ref{thm_rees}
 that there must be a generator in every $\GreenR$-class, and hence
    that $E$ contains all the non-zero idempotents of $S$; however, we
    cannot directly verify this.)

    Since $E$ is finite, we can check which words in $E$ stabilise each
    generator $a \in A$ on the right; for each $a \in A$, let $SR_a$ be the set of such.
    Check that every element of $E$ lies in $SR_b$ for some generator $b$.
    For every pair of generators $a$ and $b$, check:
    \begin{itemize}
        \item[(i)] that $SL_a$ intersects $SR_b$ in at most one element,
             and in \textbf{exactly} one element if and only if $ba$
             represents a non-zero element in the semigroup.
        \item[(ii)] that $SL_a$ and $SL_b$ are either equal or disjoint; and
        \item[(iii)] that $SR_a$ and $SR_b$ are either equal or disjoint.
    \end{itemize}
    It follows from the Rees theorem that all of these conditions must hold
    in a completely zero-simple semigroup, so if any fails, output \NO.

    Our next objective is to verify that idempotents in the same
    $SL$-class [$SR$-class] are $\GreenR$-related [$\GreenL$-related] in the semigroup.
    For this, it will suffice to verify that for every $a \in A$ and every
    pair of elements $e, f \in SL_a$ [$e,f \in SR_a$] we have $ef \GreenH f$
    [$ef \GreenH e$] in the semigroup. To verify this, we invoke
    Lemma~\ref{lemma_inverse_conditions}. If the semigroup is completely
    simple then it is easily verified that $ef$ must have finitely many
    left inverses with respect to $f$ [respectively, $e$], one of which is also right inverse.
    Hence, we can use Lemma~\ref{lemma_inverse_conditions} to check that
    $ef \GreenH f$ [respectively, $ef \GreenH e$], outputting \NO\ if it
    transpires that $ef$ has no left inverses, infinitely many left
    inverses or finitely many left inverses none of which is a right inverse,
    with respect to $f$ [respectively, $e$].

To proceed further, we employ the following lemma.

\begin{lemma}\label{lemma_findidempotents}
There is an algorithm to perform the following task:

\textbf{Instance:} a word $w \in L$ which does not represent zero;

\textbf{Problem:} \textbf{either} find an
                  idempotent in $E$ which is $\GreenL$-related to $w$ and an
                  idempotent in $E$ which is $\GreenR$-related to $w$, \textbf{or} discover that
                  the semigroup is not completely zero-simple.
\end{lemma}
\begin{proof}
First, check which idempotents in $E$ stabilise $w$ on the right. For
each such idempotent $e$, check whether there exists $q \in L$ such that
$qw = e$ in the semigroup, that is, such that $(q, e) \in L_w$. Since a
completely zero-simple semigroup has an idempotent in every $\GreenL$-class,
if there are none such, output that the semigroup is
not completely zero-simple. Otherwise,
assume we have found a non-zero idempotent $e$ which is $\GreenL$-related to $w$, and
such that $qw = e$.

Now if the semigroup is completely zero-simple then, by a standard result,
$wq$ is also a non-zero idempotent, and moreover
this idempotent must have a representative in $E$. We can locate this representative
by solving the uniform word problem; call it $f$. Now if the semigroup is
completely zero-simple then $wq = f$ is $\GreenR$-related to $w$, so that
$fw = w$. Check this condition; if it fails, output \NO. Otherwise we have
$fw = w$ and $wq = f$ so that $w \GreenR f$, as required.
\end{proof}

Using Lemma~\ref{lemma_findidempotents}, we check that every generator
$a$ is $\GreenR$-related to an idempotent (and hence to every idempotent)
in $SL_a$ and $\GreenL$-related to an idempotent in $SR_b$. In a completely
zero-simple semigroup this must be the case, so if not, output \NO.

Now, we check that for every pair of generators 
   $b$ and $a$ with non-zero product, the product $ba$ is $\GreenR$-related
   to an idempotent in $SL_b$ and $\GreenL$-related to an idempotent in $SR_a$.
   Again, we can check this by Lemma~\ref{lemma_findidempotents}, and it must
be satisfied in a completely zero-simple semigroup, so if it fails, output \NO.

Finally, we check that for every pair of generators $a,b$, there exists
a word $c_1 \dots c_n \in L$ which does not represent zero, and which has
the property that $SL_{c_1} = SL_a$ and $SR_{c_n} = SR_b$. Once again, if
this fails, output \NO.

We claim, at this point, that the semigroup is completely zero-simple.

First, we claim that for every word $a_1 \dots a_n \in A^*$, the element
represented is $\GreenR$-related to $a_1$ and $\GreenL$-related to $a_n$ in the
semigroup. If $n = 1$ there is nothing to prove. If $n = 2$ then we have
checked that $a_1 a_2$ is $\GreenL$-related to an idempotent in $SR_{a_2}$,
which in turn is $\GreenL$-related to $a_2$ itself. Similarly, $a_1 a_2$ is
$\GreenR$-related to an idempotent in $SL_{a_1}$ which in turn is
$\GreenR$-related to $a_1$ itself.

Now assume that $a_1 \dots a_n$ is a counterexample of minimal length
(necessarily 3 or more), say not $\GreenL$-related to $a_n$. Now certainly we
have $a_1 a_2 \GreenL a_2$. Now since $\GreenL$ is right compatible, we
see that $a_1 \dots a_n \GreenL a_2 \dots a_n$. But by the minimality
assumption, $a_2 \dots a_n$ is $\GreenL$-related to $a_n$. Since $\GreenL$
is transitive, this gives the required contradiction. A symmetric argument
applied in the case that $a_1 \dots a_n$ is not $\GreenR$-related to $a_1$.

It follows now that the semigroup has a single non-zero $\GreenD$-class.
Indeed, suppose $a_1 \dots a_m, b_1 \dots b_n \in A^*$ represent non-zero
elements. Then
$a_1 \dots a_m \GreenR a_1$ and $b_1 \dots b_n \GreenL b_n$. Now
there is a word $c_1 \dots c_q \in L$ which does not represent zero
such that $SL_{c_1} = SL_{a_1}$ and $SR_{c_q} = SR_{b_n}$. It follows
    that
$$a_1 \dots a_m \ \GreenR \ a_1 \ \GreenR \ c_1 \ \GreenR \ c_1 \dots c_q \ \GreenL \  c_q \ \GreenL \ b_n \ \GreenL \ b_1 \dots b_n.$$
so that $a_1 \dots a_m \GreenD b_1 \dots b_n$ as required. Thus, the semigroup
has a single non-zero $\GreenD$-class.

Moreover, we have seen that every element lies in the $\GreenR$-class
[$\GreenL$-class]
of one of the finitely many generators. So the semigroup has only
finitely many $\GreenR$-classes and $\GreenL$-classes, and hence also only finitely
many $\GreenH$-classes. Since an $\GreenH$-class can contain at most one idempotent,
we conclude that the semigroup has only finitely many idempotents, and that the
semigroup is completely zero-simple as required.
\end{proof}

A corollary is a corresponding result for completely simple semigroups.
\begin{corollary}
There is an algorithm for the following problem:

\textbf{Instance:} an interpreted automatic structure describing a
semigroup $S$;

\textbf{Problem:} decide whether $S$ is completely simple.
\end{corollary}
\begin{proof}
Clearly, $S$ is completely simple if and only if $S^0$ is completely
zero-simple. Use Proposition~\ref{prop_adjoinzero} to obtain an
automatic structure for $S^0$, and then apply Theorem~\ref{thm_compzero_decide}.
\end{proof}

\subsection{Computing the Rees matrix structure.}

In this section, we present an algorithm which, given an automatic
structure for a completely zero-simple semigroup, computes its Rees
matrix structure.

\begin{theorem}\label{thm_compzero_construction}
There is an algorithm for the following problem:

\textbf{Instance:} an interpreted automatic structure $\Gamma$ which describes
 a completely zero-simple semigroup $S$;

\textbf{Problem:} construct (i) an interpreted automatic structure $\Delta$ for the
                 maximal subgroup $G$ of $S$ and (ii) a finite 
                 matrix $Q$ with entries drawn from $L(\Delta) \cup \lbrace 0 \rbrace$, such that $S$
                 is isomorphic to the Rees matrix semigroup
                 $M^0(G; I, \Lambda; Q)$ (where $Q$ is interpreted as a
                 matrix over $G \cup \lbrace 0 \rbrace$ in the obvious way).
\end{theorem}
\begin{proof}
By Proposition~\ref{prop_uniqueness}, we may assume that the automatic
structure $\Gamma$ has uniqueness.
We begin by using the procedure from the proof of
Theorem~\ref{thm_compzero_decide} 
to find the (necessarily finite) set $E$
of representatives in $L(\Gamma)$ of idempotents. Again, we sort these
into $\GreenL$-classes and $\GreenR$-classes according to how they stabilise
each other. Let $I$ and
$\Lambda$ be the sets of non-zero $\GreenR$- and $\GreenL$-classes respectively.
For $i \in I$ and $\lambda \in \Lambda$ let $H_{i \lambda}$ denote the
$\GreenH$-class $i \cap \lambda$ and let $e_{i \lambda}$ be the word representing
the (necessarily unique) idempotent in $H_{i \lambda}$ if it exists.

Choose distinguished elements $i_0 \in I$ and $\lambda_0 \in \Lambda$
such that $H_{i_0 \lambda_0}$ contains an idempotent, that is, such
that $e_{i_0 \lambda_0}$ exists.

Next, we wish to choose for each $i \in I$ a word $r_i \in L(\Gamma)$
representing an element of $H_{i \lambda_0}$. If $H_{i \lambda_0}$
contains an idempotent then we simply set $r_i = e_{i \lambda_0}$.
Otherwise, we are involved in slightly more work. The words we seek are
exactly those of the form $awb \in L(\Gamma)$ where $a$ and $b$ are generators in
the $\GreenR$-class $i$ and the $\GreenL$-class $\lambda_0$ respectively, and $awb$
does not represent zero. Since the $\GreenH$-class cannot be empty, there must be
such a word; hence, we can find one by enumerating the set of such words
until we find one which does not represent zero. Similarly, we choose for
each $\lambda \in \Lambda$ a word $q_\lambda \in L(\Gamma)$ representing an
element of $H_{i_0 \lambda}$. Notice that by construction we have
$$r_{i_0} = q_{\lambda_0} = q_{\lambda_0} r_{i_0} = e_{\lambda_0 i_0}.$$

We now construct an $\Lambda \times I$ matrix $P$ with entries drawn
from $L(\Gamma) \cup \lbrace 0 \rbrace$, where $0$ is a new symbol. For
each $i \in I$ and $\lambda \in \Lambda$, let $P_{\lambda i}$ be the word in
$L$ representing the same element as the product
$q_\lambda r_i$. Notice that this element lies in $H_{i_0 \lambda_0}$.
It follows from the proof of \cite[Theorem~3.2.3]{Howie95} that the
semigroup $S$ is isomorphic to the Rees matrix semigroup
$M(H_{i_0 \lambda_0}; I, \Lambda; P)$, where $P$ is intepreted as
a matrix over $H_{i_0 \lambda_0} \cup \lbrace 0 \rbrace$ in the obvious
way.

It remains to construct an automatic structure for the maximal subgroup
$H_{i_0 \lambda_0}$. It follows from standard facts about completely
zero-simple semigroups that for each generator $a \in A(\Gamma)$, there exists
a unique $i_a \in I$, a unique $\lambda_a \in \Lambda$ and a unique
$b_a \in H_{\lambda_0 i_0}$ such that $a = r_{i_a} b_a s_{\lambda_a}$.
For each $a \in A(\Gamma)$, we compute (for example, by enumeration and testing,
although more efficient means are available) a representative
$w_a \in L(\Gamma)$ for $b_a$.

Define a new alphabet
$$A(\Delta) = \lbrace c_a \mid a \in A \rbrace \cup \lbrace d_{\lambda i} \mid \lambda \in \Lambda, i \in I, P_{\lambda i} \neq 0 \rbrace.$$
We view $A(\Delta)$ as an alphabet of generators for the maximal subgroup
$H_{i_0 \lambda_0}$, where each $c_a$ represents the element
$b_a \in H_{i_0 \lambda_0}$, and each $d_{\lambda i}$ represents
the element represented by the sandwich matrix entry $P_{\lambda i}$.

Let $R$ be the set of words in $L$ which represent elements of
$H_{i_0 \lambda_0}$. Define a function $\phi : R \to B^*$ by
$$\phi(a_1 a_2 \dots a_n) = c_{a_1} d_{\lambda_{a_1} i_{a_2}} c_{a_2} \dots d_{\lambda_{a_{n-1}} i_{a_n}} c_{a_n}$$
where each $i_k$ and $\lambda_k$ are such that $a_k$ represents an element
of $H_{\lambda_k i_k}$.
Let $L(\Delta) = \phi(L \cap R)$. It follows from
\cite[Proof of Theorem~4.6]{KambitesReesMatrix} that $L(\Delta)$ is a regular
language and is straightforward to compute. Now recalling that
$P_{\lambda_0 i_0} = e_{\lambda_0 i_0}$, one can show that
$$L_=(\Delta) = \lbrace (u \phi, v \phi) \mid (u, v) \in L_= \cap (R \times R) \rbrace$$
while for each $c_a \in B$ we have
$$L_{c_a}(\Delta) = \lbrace (u \phi, v \phi) \mid (u, v) \in L_{w_a} \cap (R \times R) \rbrace$$
and for each $d_{\lambda i} \in B$
$$L_{d_\lambda i}(\Delta) = \lbrace (u \phi, v \phi) \mid (u, v) \in L_{P_{\lambda i}} \cap (R \times R) \rbrace.$$
It can now be shown (for example, using arguments similar to those in
\cite[Proof of Theorem~4.6]{KambitesReesMatrix}) that each of these relations
is synchronously rational and straightforward to compute.
This completes the construction of an automatic structure for the maximal
subgroup $H_{i_0 \lambda_0}$.
 
Finally, we apply the function $\phi$ to each non-zero entry in the
sandwich matrix $P$, to obtain a new $\Lambda \times I$ matrix with
entries drawn from $K$.
\end{proof}

As an immediate corollary, we obtain a corresponding result for the case
of completely simple semigroups.

\begin{corollary}
There is an algorithm for the following problem:

\textbf{Instance:} an interpreted automatic structure $\Gamma$ which represents
 a completely simple semigroup $S$;

\textbf{Problem:} construct (i) an automatic structure $\Delta$ for the
                 maximal subgroup $G$ of $S$ and (ii) a finite 
                 matrix $Q$ with entries drawn from $L(\Delta)$, such that $S$
                 is isomorphic to the Rees matrix semigroup $M(G; I, \Lambda; Q)$ (where
                 $Q$ is interpreted as a matrix over $G$ in the obvious way).

\end{corollary}
\begin{proof}
By Proposition~\ref{prop_adjoinzero}, we can adjoin a zero to $S$ to obtain
a completely zero-simple semigroup $S^0$. By Theorem~\ref{thm_compzero_construction}
we can construct a Rees matrix representation $M^0(G; I, \Lambda; P)$ for $S^0$
together with an automatic structure for $G$. It is easily verified that
$P$ has no zero entries and $M(G; I, \Lambda; P)$ is a Rees matrix
representation for $S$.
\end{proof}

\section*{Acknowledgements}

This research was supported by a Marie Curie Intra-European Fellowship
within the 6th European Community Framework Programme. The first author
would like to thank Kirsty for all her support and encouragement.

\bibliographystyle{plain}

\end{document}